*Dedicated to Professor Roman S. Ingarden*
*for his 90th birthday*

# HYPERBOLIC HYPERCOMPLEX $D$-BAR OPERATORS, HYPERBOLIC $CR$-EQUATIONS AND HARMONICITY

## L. N. APOSTOLOVA, S. DIMIEV* and P. STOEV

*Institute of Mathematics and Informatics,*
*of the Bulgarian Academy of Sciences,*
*Sofia 1113, Bulgaria,*
*E-mail: liliana@math.bas.bg,*

*Institute of Mathematics and Informatics,*
*of the Bulgarian Academy of Sciences,*
*Sofia 1113, Bulgaria,*
*\*E-mail: sdimiev@math.bas.bg,*

*Univ. of Architecture, Civil Engineering and Geodesy,*
*Sofia , bul. Hristo Smirnenski N 1, Bulgaria*

*Univ. of Architecture, Civil Engineering and Geodesy,*
*Sofia , bul. Hristo Smirnenski N 1, Bulgaria*
*E-mail: peteruasg@abv.bg*

### Abstract

This paper is partially a review of the development of the Investigation Program announced by Stancho Dimiev at the Bedlevo Conference on Hypercomplex Analysis (2006).

A new aspect related with hyperbolic complex numbers, their generalizations and applications on hyperbolic harmonicity is included in some Appendix of S. Dimiev.

The recent results on the hyperbolic complex structures are exposed. We describe systematically the $\partial$-bar trik which led to Cauchy-Riemann type systems for different generalized hypercomplex algebras. This is in fact an elementary $\partial$-bar type theory of the PDE on the mentioned algebras. Microlocal treatment (distributions and Fourier transform) of the solutions is included too, with one result of P. Popivanov.

## 1    Introduction

The interconnection between complex analysis and partial differential equations (PDE) is well known. In the classical theory of the functions of one complex variable the notion of complex differentiability is closely related with the Cauchy-Riemann system of linear partial differential equations and the Laplacian second order operator. The progress in the theory of partial differentiable operators (PDE) related with the solution of the $\overline{\partial}$-problem of Neumann confirmed the so-called $\overline{\partial}$-technique.

Here we shortly recall the simplest version of the mentioned $\overline{\partial}$-technique. Let $f = u + iv$ be a complex value function of one complex variable $z = x + iy$, $z \in \mathbf{C}$, $x, y \in \mathbf{R}$, $u = u(x, y)$, $v = v(x, y)$, defined on on a domain $G$ in the complex plane $\mathbf{C}$. Having in mind that $z + \overline{z} = 2x$ and $z - \overline{z} = 2iy$, we introduce the following two linear differential operators with complex coefficients

$$\frac{\partial}{\partial x} - i\frac{\partial}{\partial y} \ \text{ and } \ \frac{\partial}{\partial x} + i\frac{\partial}{\partial y},$$

which act on the corresponding space of differentiable functions of two real variables $x$ and $y$. Now we denote

$$\frac{\partial}{\partial z} := \frac{1}{2}\left(\frac{\partial}{\partial x} - i\frac{\partial}{\partial y}\right) \ \text{ and } \ \frac{\partial}{\partial \overline{z}} := \frac{1}{2}\left(\frac{\partial}{\partial x} + i\frac{\partial}{\partial y}\right).$$

It is easy to see that the following identity holds for the differential $df$

$$(1) \qquad df = \frac{\partial f}{\partial x}dx + \frac{\partial f}{\partial y}dy = \frac{\partial f}{\partial z}dz + \frac{\partial f}{\partial \overline{z}}d\overline{z}.$$

The definition of complex differentiability on a domain $G \subset \mathbf{R}^2$ (or holomorphicity) of the function $f$ now seems as follows: $f$ is holomorphic on $G$ if $df = \frac{\partial f}{\partial z}dz$, or

$$(2) \qquad \frac{\partial f}{\partial \overline{z}} = \frac{\partial f}{\partial x} + i\frac{\partial f}{\partial y} = 0.$$

In view that $\frac{\partial f}{\partial x} + i\frac{\partial f}{\partial y} = \frac{\partial u}{\partial x} - \frac{\partial v}{\partial y} + i\left(\frac{\partial u}{\partial y} + \frac{\partial v}{\partial x}\right)$ we receive the classical system of Cauchy-Riemann, namely

$$(3) \qquad \frac{\partial u}{\partial x} - \frac{\partial v}{\partial y} = 0, \quad \frac{\partial u}{\partial y} + \frac{\partial v}{\partial x} = 0,$$

and the so-called Laplace equation $\Delta(U) = \frac{\partial^2 U}{\partial x^2} + \frac{\partial^2 U}{\partial y^2} = 0.$

For to continue we need some terminology:

1.) The field of complex numbers considered as a complex division algebra is always denoted by $\mathbf{C}$.

2.) There are different names for the corresponding non-division algebra received after changing $i^2 = -1$ by $j^2 = 1$, i. e. $x + jy$, $j^2 = 1$, $j \notin \mathbf{C}$.

a) algebra of double-real numbers (Rosenfeld)

b) paracomplex algebra according P. Lieberman

c) hyperbolic complex numbers algebra

We shall use the third term (the case c), denoting the mentioned algebra by $\tilde{\mathbf{C}}$.

So, $\mathbf{C}$ is the complex division algebra of ordinary complex numbers ($i^2 = -1$) and $\tilde{\mathbf{C}}$ is the corresponding non-division algebra of hyperbolic complex numbers ($j^2 = 1$, $j \notin \mathbf{C}$).

We shall consider functions $f : G \to \tilde{\mathbf{C}}$, $G \subset \tilde{\mathbf{C}}$ or $f(z) = f_0(x, y) + jf_1(x, y)$, $z \in G$, $z = x + jy$, $j^2 = 1$, $f_0, f_1 : G \to \mathbf{R}$. In the next we always suppose that the considered functions are continuously differentiables. The following representation for the differential $df$ is valid

$$(4) \qquad df = \frac{\partial f}{\partial x}dx + \frac{\partial f}{\partial y}dy = \frac{\partial f}{\partial z}dz + \frac{\partial f}{\partial z^*}dz^*,$$

where $\frac{\partial}{\partial z} := \frac{1}{2}\left(\frac{\partial}{\partial x} + j\frac{\partial}{\partial y}\right)$, $\frac{\partial}{\partial z^*} := \frac{1}{2}\left(\frac{\partial}{\partial x} - j\frac{\partial}{\partial y}\right)$, $\frac{1}{j} = j$, $z^* := x - jy$, and $dz = dx + jdy$, $dz^* = dx - jdy$.

By definition, the function $f$ is called differentiable on $G$ if

$$(5) \qquad \frac{\partial f}{\partial z^*} = \frac{\partial f}{\partial x} - j\frac{\partial f}{\partial y} = 0.$$

The condition above implies that

$$\frac{\partial f_0}{\partial x} - \frac{\partial f_1}{\partial y} + j\left(\frac{\partial f_1}{\partial x} - \frac{\partial f_0}{\partial y}\right) = 0.$$

Finally, (as $j^2 = 1$) we get the following linear system of PDE

$$(6) \qquad \frac{\partial f_0}{\partial x} - \frac{\partial f_1}{\partial y} = 0, \qquad \frac{\partial f_0}{\partial y} - \frac{\partial f_1}{\partial x} = 0.$$

This is an analogous of the Cauchy-Riemann system. Respectively, we receive the hyperbolic second order PDE

$$(7) \qquad \frac{\partial^2 U}{\partial x^2} - \frac{\partial^2 U}{\partial y^2} = 0,$$

satisfied by $f_0$ and $f_1$.

Using the geometric interpretation of the above considered algebras as 2-dimensional planes (as real ones), we make the summary:

on the complex plane $\mathbf{C}$ the complex differentiability on a domain in $\mathbf{C}$ is expressed by the linear (elliptic) first order Cauchy-Riemann system and the second order Laplace operator $\Delta$, on the hyperbolic plane of the hyperbolic numbers the analogous differentiation is expressed by a linear hyperbolic system of first order and a hyperbolic equation of second order.

All this is well known. Our purpose is to develop the same question for generalized hypercomplex algebras (fourth-real, bicomplex, coquaternion's etc.)

## 2 Fourth-real functions of hyperbolic cyclic 4-real variables [DMZ]

The algebra $\mathbf{R}(1, j, j^2, j^3)$ of the hyperbolic cyclic 4-real variables $\zeta = x + yj + uj^2 + vj^3$, $x, y, u, v \in \mathbf{R}$, $j^4 = 1$, is a natural generalization of the algebra of hyperbolic 2-real (hyperbolic complex) variables $R(1, j)$, $j^2 = 1$. Each hyperbolic 4-real function $\Phi : R(1, j, j^2, j^3) \rightarrow R(1, j, j^2, j^3)$ admits a representation as a function of two hyperbolic 2-real variables, namely $\Phi(\zeta) = \Phi_0(\alpha, \beta) + j\Phi_1(\alpha, \beta)$, where $\zeta = \alpha + j\beta$, $\alpha = x + j^2 u$, $\beta = y + j^2 v$, $\alpha^* = x - j^2 u$, $\beta^* = y - j^2 v$.

Let $f$ be a 4-real function of 4-real variables, $f = f_0 + f_1 + f_2 j^2 + f_3 j^3$, the coordinate functions $f_k (k = 0, 1, 2, 3)$ being real-valued functions of 4 real variables. Having the differential $df = \dfrac{\partial f}{\partial x} dx + \dfrac{\partial f}{\partial y} dy + \dfrac{\partial f}{\partial u} du + \dfrac{\partial f}{\partial v} dv$ we introduce the following differential operators $\dfrac{\partial}{\partial \alpha}, \dfrac{\partial}{\partial \beta}, \dfrac{\partial}{\partial \alpha^*}, \dfrac{\partial}{\partial \beta^*}$ :

$$\frac{\partial f}{\partial \alpha} := \frac{1}{2}\left(\frac{\partial f}{\partial x} + j^2 \frac{\partial f}{\partial u}\right), \quad \frac{\partial}{\partial \beta} := \frac{1}{2}\left(\frac{\partial f}{\partial y} + j^2 \frac{\partial f}{\partial v}\right),$$

$$\frac{\partial f}{\partial \alpha^*} := \frac{1}{2}\left(\frac{\partial f}{\partial x} - j^2 \frac{\partial f}{\partial u}\right), \quad \frac{\partial}{\partial \beta^*} := \frac{1}{2}\left(\frac{\partial f}{\partial y} - j^2 \frac{\partial f}{\partial v}\right).$$

One calculate that the differential $df$ above will seem as follows

(4′)
$$df = \frac{\partial f}{\partial \alpha} d\alpha + \frac{\partial f}{\partial \beta} d\beta + \frac{\partial f}{\partial \alpha^*} d\alpha^* + \frac{\partial f}{\partial \beta^*} d\beta^*.$$

One calculate also that

$$\frac{\partial f}{\partial \alpha^*} d\alpha^* + \frac{\partial f}{\partial \beta^*} d\beta^* = \frac{1}{2}\left(\frac{\partial f}{\partial u} - j^2 \frac{\partial f}{\partial x}\right) du + \frac{1}{2}\left(\frac{\partial f}{\partial v} - j^2 \frac{\partial f}{\partial y}\right) dv.$$

The $\partial$-bar trick is defined here by the equation

(5′)
$$\frac{\partial f}{\partial \alpha^*} d\alpha^* + \frac{\partial f}{\partial \beta^*} d\beta^* = 0.$$

As a result we obtain that the function $f$ is a four-real differentiable on its domain of definition if the following two equations hold

(6′)
$$\frac{\partial f}{\partial x} = j^2 \frac{\partial f}{\partial u}, \quad \frac{\partial f}{\partial y} = j^2 \frac{\partial f}{\partial v}, \quad j^4 = +1.$$

In terms of coordinates functions $f_0, f_1, f_2, f_3$ we receive a system of 8 first order hyperbolic partial differential equations (see [DMJ]). They are the following ones

$$\frac{\partial f_0}{\partial x} = \frac{\partial f_2}{\partial u}, \qquad \frac{\partial f_0}{\partial u} = \frac{\partial f_2}{\partial x}, \qquad \frac{\partial f_1}{\partial x} = \frac{\partial f_3}{\partial u}, \qquad \frac{\partial f_1}{\partial u} = \frac{\partial f_3}{\partial x},$$

$$\frac{\partial f_0}{\partial y} = \frac{\partial f_2}{\partial v}, \qquad \frac{\partial f_0}{\partial v} = \frac{\partial f_2}{\partial y}, \qquad \frac{\partial f_1}{\partial y} = \frac{\partial f_3}{\partial v}, \qquad \frac{\partial f_1}{\partial v} = \frac{\partial f_3}{\partial y}.$$

## 3 Bicomplex and double-complex algebras

These two algebras are in fact isomorphic. The first one, denoted here by $BC$, was developed in the book [BP]. Its elements are the bicomplex numbers $\alpha = z + jw$, with $j^2 = -1$, and the elements of the second are of the same $\alpha = z + jw$, but with $j^2 = i$, $i \in \mathbf{C}$. A motivation for the second, denoted by $DC$, is related with the corresponding function theory. In the case of bicomplex functions $f : BC \to BC$ one present $f$ in the form $f(\alpha) = f_0(z) + jf_1(w)$, [D], [RS] and others. In the second one, the functions $f : DC \to DC$ are presented in the form $f(\alpha) = f_0(z, w) + jf_1(z, w)$, which is more natural in view of the possibility to apply the theory of the functions of two complex variables directly to $f_0(z, w)$ and $f_1(z, w)$.

The $\partial$-bar trick for bicomplex functions is easily obtained with the help of the transformation

$$\frac{\partial f}{\partial \alpha} = \frac{1}{2}\left(\frac{\partial f}{\partial z} - j\frac{\partial f}{\partial w}\right), \quad \frac{\partial f}{\partial \alpha^*} = \frac{1}{2}\left(\frac{\partial f}{\partial z} + j\frac{\partial f}{\partial w}\right), \quad \alpha^* = z - jw.$$

Indeed, we calculate that $\dfrac{\partial f}{\partial \alpha}d\alpha + \dfrac{\partial f}{\partial \alpha^*}d\alpha^* = \dfrac{\partial f}{\partial z}dz - j^2\dfrac{\partial f}{\partial w}dw$, which gives for $j^2 = -1$, just the differential $df = \dfrac{\partial f}{\partial z}dz + \dfrac{\partial f}{\partial w}dw$, or the identity

$$(4'') \qquad df = \frac{\partial f}{\partial z}dz + \frac{\partial f}{\partial w}dw = \frac{\partial f}{\partial \alpha}d\alpha + \frac{\partial f}{\partial \alpha^*}d\alpha^*.$$

So, for to have $df = \dfrac{\partial f}{\partial \alpha}d\alpha$ we must set

$$(5'') \qquad \frac{\partial f}{\partial \alpha^*} = \frac{1}{2}\left(\frac{\partial f}{\partial z} + j\frac{\partial f}{\partial w}\right) = 0.$$

The bicomplex Cauchy-Riemann system for $f_0$ and $f_1$ follows immediately:

$$(6'') \qquad \frac{\partial f_0}{\partial z} = \frac{\partial f_1}{\partial w}, \quad \frac{\partial f_0}{\partial w} = -\frac{\partial f_1}{\partial z}.$$

For the double-complex algebra $DC$ the above used $\dfrac{\partial f}{\partial \alpha}$, $\dfrac{\partial f}{\partial \alpha^*}$ are slightly different

$$\frac{\partial f}{\partial \alpha} = \frac{1}{2}\left(\frac{\partial f}{\partial z} - ij\frac{\partial f}{\partial w}\right), \quad \frac{\partial f}{\partial \alpha^*} = \frac{1}{2}\left(\frac{\partial f}{\partial z} + ij\frac{\partial f}{\partial w}\right).$$

The holomorphicity condition is the same

$$(5''') \qquad \frac{\partial f}{\partial \alpha^*} = \frac{1}{2}\left(\frac{\partial f}{\partial z} + ij\frac{\partial f}{\partial w}\right) = 0,$$

and the corresponding double-complex Cauchy-Riemann system seems as follows

$$(6''') \qquad \frac{\partial f_0}{\partial z} = \frac{\partial f_1}{\partial w}, \quad \frac{\partial f_0}{\partial w} = i\frac{\partial f_1}{\partial z}.$$

The second order partial differential equations

$$\frac{\partial^2 U}{\partial z^2} + \frac{\partial^2 U}{\partial w^2} = 0 \ \text{ and } \ \frac{\partial^2 U}{\partial z^2} + i\frac{\partial^2 U}{\partial w^2} = 0,$$

are satisfied by $f_0$ and $f_1$ in the bicomplex case (the first one - the bicomplex Laplacian) and respectively in the second case (the second one - the double-complex Laplacian). The mentioned "Laplacians" are partially hypoellipic (to see the next).

# 4 Hyperbolic bicomplex numbers

The algebraic property of bicomplex numbers and hyperbolic numbers have been study in [RS], [Sg] and [Y]. In this note we develop the notion of hyperbolic bicomplex numbers (variables) and the basic notion of the hyperbolic bicomplex analysis, especially hyperbolic complex CR system and the corresponding "hyperbolic" complex Laplacian (see also [AKK]).

## 4.1 Definition of the hyperbolic bicomplex numbers

First we recall the definition of the bicomplex numbers according to B. Price. Let $\alpha = z + j_2 w$ be a bicomplex number, i.e. $z, w \in \mathbf{C}$, $z = x_0 + j_1 x_1$, $w = x_2 + j_1 x_3$, $x_0, x_1, x_2, x_3 \in \mathbf{R}$, $j_1^2 = j_2^2 = -1$ and $j_1 j_2 = j_2 j_1$.

We say that $\alpha = z + j_2 w$ is a hyperbolic bicomplex number if

$$(11) \qquad j_1^2 = j_2^2 = +1 \text{ and } j_1 j_2 = j_2 j_1.$$

The above given definition implies that $z$ and $w$ become hyperbolic complex numbers, sometimes denoted $\tilde{z}$ and $\tilde{w}$.

The conjugate of the bicomplex number $\alpha$, denoted by $\alpha^*$, is by definition

$$(12) \qquad \alpha^* = z - j_2 w.$$

Analogously, the conjugate of the hyperbolic bicomplex number $\tilde{\alpha}$, denoted $\tilde{\alpha}^*$, is by definition

$$(13) \qquad \tilde{\alpha}^* = \tilde{z} - j_2 \tilde{w}.$$

Clearly, we have

$$(14) \qquad \tilde{\alpha} \tilde{\alpha}^* = \tilde{z}^2 - \tilde{w}^2.$$

The set of all bicomplex numbers is a commutative ring with zero divisors. The same is true for the hyperbolic bicomplex numbers.

**Lemma.** *The hyperbolic bicomplex numbers $\alpha = \zeta + \eta j_2$ is a zero divisor iff*

$$(15) \qquad \zeta^2 - \eta^2 = 0.$$

P r o o f. We have the identity

$$(\zeta + j_2 \eta)(\zeta - j_2 \eta) = \zeta^2 - j_2^2 \eta^2 = \zeta^2 - \eta^2. \qquad \Diamond$$

## 4.2 Hyperbolic bicomplex functions of hyperbolic bicomplex variable

Following the notations of B. Price $\mathbf{C_2} = \{z + jw : z, w \in \mathbf{C}, j^2 = -1\}$, we denote the hyperbolic bicomplex functions as follows

$$f : \tilde{\mathbf{C_2}} \to \tilde{\mathbf{C_2}}, \ \tilde{\mathbf{C_2}} = \{\tilde{\zeta} + \tilde{\eta} j_2 : \tilde{\zeta} \in \tilde{\mathbf{C}}, \tilde{\eta} \in \tilde{\mathbf{C}}, j_2^2 = +1\},$$

where $\tilde{\mathbf{C}}_2$ is the algebra of the hyperbolic bicomplex numbers. In fact we have

$$(16) \qquad f(\tilde{\alpha}) = f_0(\tilde{\zeta}, \tilde{\eta}) + f_1(\tilde{\zeta}, \tilde{\eta})j_2, \quad \tilde{\alpha} = \tilde{\zeta} + j_2\tilde{\eta}.$$

where $\tilde{\zeta} = x_0 + j_1x_1$, $\tilde{\eta} = x_2 + j_1x_3$, $j_1^2 = +1$. Changing the notations as follows $z := \tilde{\zeta}$ and $w := \tilde{\eta}$ we have

$$(17) \qquad f(\tilde{\alpha}) = f_0(z, w) + f_1(z, w)j_2,$$

with $z = x_0 + j_1x_1$ and $w = x_2 + j_1x_3$ with $j_1^2 = +1$.

Let us take the differential

$$(18) \qquad df = \frac{\partial f}{\partial z}dz + \frac{\partial f}{\partial w}dw,$$

where

$$\frac{\partial f}{\partial z} := \frac{\partial f_0}{\partial z} + \frac{\partial f_1}{\partial z}j_2, \qquad \frac{\partial f}{\partial w} := \frac{\partial f_0}{\partial w} + \frac{\partial f_1}{\partial w}j_2.$$

**Proposition.** *The following identity is valid*

$$(19) \qquad \frac{\partial f}{\partial z}dz + \frac{\partial f}{\partial w}dw = \frac{\partial f}{\partial \tilde{\alpha}}d\tilde{\alpha} + \frac{\partial f}{\partial \tilde{\alpha}^*}d\tilde{\alpha}^*,$$

*where*

$$\frac{\partial f}{\partial \tilde{\alpha}} := \frac{1}{2}\left(\frac{\partial f}{\partial z} + j_2\frac{\partial f}{\partial w}\right), \qquad \frac{\partial f}{\partial \tilde{\alpha}^*} := \frac{1}{2}\left(\frac{\partial f}{\partial z} - j_2\frac{\partial f}{\partial w}\right),$$

$$d\tilde{\alpha} = dz + j_2dw, \qquad d\tilde{\alpha}^* = dz - j_2dw.$$

P r o o f. By direct calculation. $\qquad\qquad\qquad\qquad\qquad\qquad \diamond$

We say that the function $f$ is a hyperbolic holomorphic bicomplex function if

$$(20) \qquad \frac{\partial f}{\partial \tilde{\alpha}^*} = 0.$$

**Theorem.** *The function $f$ is holomorphic hyperbolic bicomplex one iff the system*

$$(21) \qquad \frac{\partial f_0}{\partial z} = \frac{\partial f_1}{\partial w}, \qquad \frac{\partial f_0}{\partial w} = \frac{\partial f_1}{\partial z},$$

*is satisfied by $f_0$ and $f_1$.*

P r o o f. The condition (20) implies

$$\frac{\partial f}{\partial z} = \frac{\partial f}{\partial w}j_2,$$

or

$$\frac{\partial f_0}{\partial z} + \frac{\partial f_1}{\partial z}j_2 = \left(\frac{\partial f_0}{\partial w}j_2 + \frac{\partial f_1}{\partial w}\right).$$

Finally, we get the system (21). $\qquad\qquad\qquad\qquad\qquad\qquad\qquad \diamond$

**Remark.** The system (21) is called hyperbolic Cauchy-Riemann bicomplex system. Let us recall that the (elliptic) Cauchy-Riemann bicomplex system seems as follows (B. Price)

(22)
$$\frac{\partial f_0}{\partial z} = \frac{\partial f_1}{\partial w}, \qquad \frac{\partial f_0}{\partial w} = -\frac{\partial f_1}{\partial z}.$$

Respectively, we obtain that for the holomorphic hyperbolic bicomplex function $f = f_0 + f_1 j_2$ the following partial differential equation of second order is satisfied by its even and odd parts

$$\frac{\partial^2 f_0}{\partial z^2} - \frac{\partial^2 f_0}{\partial w^2} = 0 \ \text{ and } \ \frac{\partial^2 f_1}{\partial z^2} - \frac{\partial^2 f_1}{\partial w^2} = 0.$$

We call the operator

$$\Delta^h = \frac{\partial^2}{\partial z^2} - \frac{\partial^2}{\partial w^2}$$

"hyperbolic" Laplace operator.

## 5  Hyperbolic double-complex functions

Here we expose once again one paragraph of the paper [ADMS] in order to correct the bad formulation of the definition of the operators $\dfrac{\partial}{\partial q}$ and $\dfrac{\partial}{\partial q^*}$.

We denote by $D\tilde{\mathbf{C}}$ the algebra of hyperbolic double-complex matrix numbers (variables) $D\tilde{\mathbf{C}} = \{Z + jW : Z, W \in \tilde{\mathbf{C}}, \ j^2 = \mathbf{j}, \ \mathbf{j}^2 = 1\}$. Let $q = \zeta + j\eta$ be a hyperbolic double-complex matrix number. Then $f(q) = f_0(\zeta, \eta) + f_1(\zeta, \eta)j$.

**Proposition.** *The next formula is valid* $\dfrac{\partial f}{\partial q}dq + \dfrac{\partial f}{\partial q^*}dq^* = \dfrac{\partial f}{\partial \zeta}d\zeta + \dfrac{\partial f}{\partial \eta}d\eta$, *where*

$$\frac{\partial f}{\partial q} = \frac{1}{2}\left(\frac{\partial f}{\partial \zeta} + \frac{1}{j}\frac{\partial f}{\partial \eta}\right), \ \ \frac{\partial f}{\partial q^*} = \frac{1}{2}\left(\frac{\partial f}{\partial \zeta} - \frac{1}{j}\frac{\partial f}{\partial \eta}\right),$$

$$dq = d\zeta_0 + jd\eta, \ \ dq^* = d\zeta - jd\eta \ .$$

P r o o f. It is enough to make addition of the following two expressions

$$\frac{\partial f}{\partial q}dq = \frac{1}{2}\frac{\partial f}{\partial \zeta}d\zeta + \frac{1}{2}\frac{\partial f}{\partial \eta}d\eta + \frac{1}{2}\frac{\partial f}{\partial \zeta}jd\eta + \frac{1}{2j}\frac{\partial f}{\partial \eta}jd\zeta,$$

$$\frac{\partial f}{\partial q^*}dq^* = \frac{1}{2}\frac{\partial f}{\partial \zeta}d\zeta + \frac{1}{2}\frac{\partial f}{\partial \eta}d\eta - \frac{1}{2}\frac{\partial f}{\partial \zeta}jd\eta - \frac{1}{2j}\frac{\partial f}{\partial \eta}jd\eta. \qquad \Diamond$$

**Definition.** We say that the function $f$ is holomorphic if and only if

$$\frac{\partial f}{\partial q^*} = \frac{1}{2}\left(\frac{\partial f}{\partial \zeta} - \frac{1}{j}\frac{\partial f}{\partial \eta}\right) = 0.$$

**Theorem.** *The function $f$ is holomorphic hyperbolic double-complex function if and only if the system* $\dfrac{\partial f_0}{\partial \zeta} = \dfrac{\partial f_1}{\partial \eta}, \quad \dfrac{\partial f_0}{\partial \eta} = \mathbf{j}\dfrac{\partial f_1}{\partial \zeta}$ *is satisfied by $f_0$ and $f_1$.*

P r o o f .  By the holomorphicity condition it follows $\dfrac{\partial f}{\partial \zeta} = \dfrac{1}{j}\dfrac{\partial f}{\partial \eta}$. But $f = f_0 + jf_1$, from where

$$\frac{\partial f_0}{\partial \zeta} + j\frac{\partial f_1}{\partial \zeta} = j^3\left(\frac{\partial f_0}{\partial \eta} + j\frac{\partial f_1}{\partial \eta}\right), \quad \text{or} \quad \frac{\partial f_0}{\partial \zeta} + j\frac{\partial f_1}{\partial \zeta} = j\mathbf{j}\frac{\partial f_0}{\partial \eta} + \frac{\partial f_1}{\partial \eta}.$$

After comparison with respect to $j$ we received the announcement of the theorem. $\diamond$

## 6  Coquaternionic complex type functions

The coquaternion algebra is an associative non-commutative real algebra with zero divisors introduced by James Cookle (1849). A coquaternion $q$ is defined as

$$q = x_0 + x_1 i + x_2 j + x_3 k, \ x_s \in \mathbf{R}, \ x_s = 0, 1, 2, 3,$$

where $i^2 = -1, j^2 = k^2 = +1, ij = k = -ij, j \notin \mathbf{C}, k \notin \mathbf{C}$.

We shall consider coquatenion congugate $q^*$ of the following form

$$q^* = x_0 + x_1 i - (x_2 + x_3 i)j,$$

different from the ordinary quaternin (1,3)-conjugate. Setting $z := x_0 + x_1 i$, $w := x_2 + x_3 i$ we receive $q = z + wj$, $q^* = z - wj$, the product $qq^*$ being a not positively defined (2,2) form. In this kind of complex representation the considered coquaternions are called complex type coquaternions (in right complex representation type) and its algebra, denoted $coHC$, is called coquaternion complex type algebra. Let us remark that coquaternions are known also under the term split-quaternions and also pseudoquaternions [Ros].

We introduce also the complex type $\partial$-bar conjugation

$$q = z + wj \to \overline{q} := \overline{z} + \overline{w}j$$

and the corresponding algebra will be denoted $coH\overline{C}$. We need of the cartesian product $coHC \times coH\overline{C} = \{q, \overline{q}\} = \{z + wj, \overline{z} + \overline{w}j\}$.

We shall consider smooth functions $f : coHC \times coH\overline{C} \to coHC \times coH\overline{C}$ or

$$f = f(q, \overline{q}) = f_0(z, w, \overline{z}, \overline{w}) + f_1(z, w, \overline{z}, \overline{w})j, \quad j^2 = +1,$$

where $f_0$ and $f_1$ are smooth functions, which we assume to be continuously partial differentiable. Its differentials are considered as a differential of smooth functions of complex variables

$$df_0 = \frac{\partial f_0}{\partial z}dz + \frac{\partial f_0}{\partial w}dw + \frac{\partial f_0}{\partial \overline{z}}d\overline{z} + \frac{\partial f_0}{\partial \overline{w}}d\overline{w},$$

$$df_1 = \frac{\partial f_1}{\partial z}dz + \frac{\partial f_1}{\partial w}dw + \frac{\partial f_1}{\partial \overline{z}}d\overline{z} + \frac{\partial f_1}{\partial \overline{w}}d\overline{w}.$$

So, we obtain $df(q, \overline{q}) = df_0(z, w, \overline{z}, \overline{w}) + df_0(z, w, \overline{z}, \overline{w})j$, or $df = df_0 + df_1 j$.

We need of the following new kind of differentials. The following 1-form of complex variables $z, \overline{z}, w, \overline{w}$ is called a non-commutative differential with respect to the unity $j$:

$$d^{(j)}f(q,\overline{q}) = \frac{\partial f}{\partial z}dz + \frac{\partial f}{\partial \overline{z}}d\overline{z} + \frac{\partial f}{\partial w}(jdwj) + \frac{\partial f}{\partial \overline{w}}(jd\overline{w}j), \tag{23}$$

where $\frac{\partial f}{\partial z} = \frac{\partial f_0}{\partial z} + \frac{\partial f_1}{\partial z}j$, $\frac{\partial f}{\partial \overline{z}} = \frac{\partial f_0}{\partial \overline{z}} + \frac{\partial f_1}{\partial \overline{z}}j$ etc. for $w$ and $\overline{w}$.

Let us remark that in the case $j$ to be replace by the basic unit 1 the ordinary notion of differential is received.

It is easy to calculate that the non-commutative differential (23) (in $dz$, $d\overline{z}$, $dw$, $d\overline{w}$) can be represented as follows

$$d^{(j)}f(q,\overline{q}) = \frac{\partial f}{\partial q}dq + \frac{\partial f}{\partial \overline{q}}d\overline{q} + \frac{\partial f}{\partial q^*}dq^* + \frac{\partial f}{\partial \overline{q}^*}d\overline{q}^*, \tag{24}$$

where

$$\frac{\partial f}{\partial q} = \frac{1}{2}\left(\frac{\partial f}{\partial z} + \frac{\partial f}{\partial w}j\right), \quad \frac{\partial f}{\partial \overline{q}} = \frac{1}{2}\left(\frac{\partial f}{\partial \overline{z}} + \frac{\partial f}{\partial \overline{w}}j\right),$$

$$\frac{\partial f}{\partial q^*} = \frac{1}{2}\left(\frac{\partial f}{\partial z} - \frac{\partial f}{\partial w}j\right), \quad \frac{\partial f}{\partial \overline{q}^*} = \frac{1}{2}\left(\frac{\partial f}{\partial \overline{z}} - \frac{\partial f}{\partial \overline{w}}j\right),$$

$$dq = dz + dw, \quad d\overline{q} = d\overline{z} + d\overline{w}, \quad dq^* = dz - dw, \quad d\overline{q}^* = d\overline{z} - d\overline{w}.$$

The smooth, respectively the complex holomorphic, coquaternion complex type function $f$ is called (*)-regular if its non-commutative differential does not contains members with $\overline{q}$ and $\overline{q}^*$, i.e.

$$d^{(j)}f(q,\overline{q}) = \frac{\partial f}{\partial q}dq + \frac{\partial f}{\partial q^*}dq^*. \tag{25}$$

This means that $\partial$-bar trick is formulated here as

$$\frac{\partial f}{\partial \overline{q}}d\overline{q} + \frac{\partial f}{\partial \overline{q}^*}d\overline{q}^* = 0. \tag{26}$$

The following theorem is valid [DKT]

**Theorem.** *The even part $f_0$ and the odd part $f_1$ of the coquaternion complex type smooth function $f$ satisfy the following system of complex partial differential equations*

$$\frac{\partial f_0}{\partial z} = \frac{\partial f_1}{\partial w}, \quad \frac{\partial f_0}{\partial w} = \frac{\partial f_1}{\partial z}, \tag{27}$$

$$\frac{\partial f_0}{\partial \overline{z}} = \frac{\partial f_1}{\partial \overline{w}}, \quad \frac{\partial f_0}{\partial \overline{w}} = \frac{\partial f_1}{\partial \overline{z}}. \tag{$\overline{27}$}$$

*For holomorphic functions $f_0$ and $f_1$ the above written system (27) and ($\overline{27}$) reduces to the unique system (27), as the second one ($\overline{27}$) is trivially satisfied.*

For the proof to see [DKT].

**Corollary.** In the smooth case the following second order system is valid

$$\frac{\partial^2 U}{\partial z^2} - \frac{\partial^2 U}{\partial w^2} = 0, \quad \frac{\partial^2 U}{\partial \overline{z}^2} - \frac{\partial^2 U}{\partial \overline{w}^2} = 0. \tag{28}$$

For holomorphic functions $f_0$ and $f_1$ the above written system reduces to its first component, which appear as a kind of "hyperbolic" second order PDE.

# 7 Fundamental solution of the hyperbolic complex $\partial$-bar operator

Concrete examples of hyperholomophic $\partial$-bar operators was obtained above for different generalized hyper-complex algebras. Here we shall study the fundamental solutions of some of them. In this paragraph we shall use elements of the theory of distributions. We recall here that a distribution is a continuous functional on the so-called test function space which elements are denoted by $\varphi$. Distributions are considered as extension of ordinary differentiable functions, but every distribution is differentiable and the operation of the differentiation is a continuous operation on the space of distributions, denoted by $\mathcal{D}'$, where$\mathcal{D}$ is the test function space. We will consider a space of test functions on $\mathbf{R}^n$ called Schwartz space that consists of smooth, rapidly decreasing functions. In the case of double-complex Laplacian $\Delta_+$ we have in mind $\mathbf{R}^4$.

In the next by $D^\infty(G)$ is denoted the space of $C^\infty$-differentiable functions defined on a domain $G$. Let $P(D)$ denotes a partially differential operator which acts on $D$.

We recall that a distribution $\varepsilon(x) \in \mathcal{D}'$ is a fundamental solution for $P(D)$ if

$$P(D)\varepsilon(x) = \delta(x), \ x \in \mathbf{R}^4, \ \delta(x) \ \text{being the Dirac delta function.}$$

Let us consider the case of the algebra of hyperbolic 2-real numbers $\mathbf{R}(1, j)$, $j^2 = +1$, $j \neq \pm1$, i.e. the algebra of hyperbolic complex numbers. We shall search a fundamental solution $E$ for the hyperbolic $\partial$-bar operator $\partial^*$

$$\partial^* = \frac{\partial}{\partial z^*} = \frac{1}{2}\left(\frac{\partial}{\partial x} - j\frac{\partial}{\partial y}\right), \quad z = x + jy, \quad x, y \in \mathbf{R}, \ j^2 = +1.$$

We change the coordinate system $(x, y)$ by a hyperbolic one $(r, t)$ as follows

$$x = r\cosh t, \quad y = r\sinh t$$

where $r > 0, r, t \in \mathbf{R}$.

It is fulfilled the rule of differentiation of compose function, by where

$$\frac{\partial}{\partial r} = \frac{\partial}{\partial x}\frac{\partial x}{\partial r} + \frac{\partial}{\partial y}\frac{\partial y}{\partial r} = \cosh t\frac{\partial}{\partial x} + \sinh t\frac{\partial}{\partial y},$$

$$\frac{\partial}{\partial t} = \frac{\partial}{\partial x}\frac{\partial x}{\partial t} + \frac{\partial}{\partial y}\frac{\partial y}{\partial t} = r\sinh t\frac{\partial}{\partial x} + r\cosh t\frac{\partial}{\partial y}.$$

Solving these equalities with respect to $\dfrac{\partial}{\partial x}$ and $\dfrac{\partial}{\partial y}$ we obtain the equalities

$$\frac{\partial}{\partial x} = \cosh t\frac{\partial}{\partial r} - \frac{1}{r}\sinh t\frac{\partial}{\partial t},$$

$$\frac{\partial}{\partial y} = -\sinh t\frac{\partial}{\partial r} + \frac{1}{r}\cosh t\frac{\partial}{\partial t}.$$

Computing the formal derivative $\frac{1}{2}\left(\frac{\partial}{\partial x} - j\frac{\partial}{\partial y}\right)$ we obtain

$$\frac{1}{2}\left(\frac{\partial}{\partial x} - j\frac{\partial}{\partial y}\right) = \frac{1}{2}e^{jt}\left(\frac{\partial}{\partial r} - \frac{j}{r}\frac{\partial}{\partial t}\right),$$

where the exponential hyperbolic function $e^z$, defined by the power series $\sum_{k=0}^{\infty}\frac{z^k}{k!}$, satisfies the Euler formula $e^{jt} = \cosh t + j\sinh t$.

A distribution $E(z)$ is called a fundamental solution of the partial differential equation $Lu = 0$ with hyperbolic complex variables, if it is fulfilled $LE(z) = \delta(z)$, where $\delta(z)$ is the Dirac delta function on $\tilde{\mathbf{C}} = \mathbf{R}(1, j)$.

**Theorem.** *The function $E(z) = \dfrac{2}{z}$ is a fundamental solution of the hyperbolic complex $\partial$-bar operator $\dfrac{\partial}{\partial z^*} = \dfrac{1}{2}\left(\dfrac{\partial}{\partial x} - j\dfrac{\partial}{\partial y}\right)$.*

**Remark.** The number $z^* = x - jy$ is a hyperbolic complex number, which is called a conjugate hyperbolic complex number for the hyperbolic complex number $z = x + jy$.

P r o o f: We follow the proof exposed in the book of L. Schwartz, Téorie des distribution, t. I, Paris, 1957, p. 48.

First, let us note that the function $E(z)$ is a locally integrable with respect to the Lesbeages measure $dxdy$ function on $\mathbf{R}^2$, hyperbolic holomorphic on the plane $\tilde{\mathbf{C}}$ except the origin, i.e. $E(z)$ satisfy the equation $\frac{\partial E(z)}{\partial z^*} = 0$ for $z \neq 0$.

Then we compute the integral on the unbounded domain

$$D := \{(r\cosh t, r\sinh t) : 0 < r < \infty, \ -\infty < t < \infty, \ r, \ t \in \mathbf{R}\}$$

$$-\left\langle E, \frac{\partial\varphi}{\partial z^*}\right\rangle = -\int\int_D \frac{2}{x + jy}\frac{\partial\varphi}{\partial z^*}\,dxdy.$$

Let us change the cartesian coordinates $x, y$ in $\mathbf{R}^2$ with the hyperbolic one above. It is fulfilled $dxdy = rdrdt$ and $z^{-1} = r^{-1}e^{-jt}$. Using the computation of the $\partial^*$-operator in coordinates $(r, t)$ we obtain

$$-\left\langle E, \frac{\partial\varphi}{\partial z^*}\right\rangle = \frac{2}{2}\int_0^\infty dr\int_{-\infty}^\infty \left(\frac{\partial}{\partial r}\cosh t - \frac{j}{r}\frac{\partial}{\partial t}\right)\varphi(r\cosh t, r\sinh t)\,dt.$$

It is true that for every infinitely smooth rapidly decreasing in infinity function $u(t)$ on the real axis $-\infty < t < \infty$, is fulfilled $\int_{-\infty}^\infty \frac{\partial u(t)}{\partial t}\,dt = 0$.

Let us denote by $\varphi^\natural(r)$ the mean value of the function $\varphi$ on the hyperbola $(r\cosh t, r\sinh t)$ with $t \in (-\infty, \infty)$ for $r \in (0, 1]$ i. e.

$$\varphi^\natural(r) = \int_{-\infty}^\infty \varphi(r\cosh t, r\sinh t)\,dt.$$

As the equality $\left(\frac{\partial\varphi}{\partial r}\right)^\natural = \frac{\partial\varphi^\natural}{\partial r}$ is fulfilled we obtain

$$-\left\langle E, \frac{\partial\varphi}{\partial z^*}\right\rangle = -\lim_{\varepsilon\to 0}\int_\varepsilon^\infty \frac{\partial\varphi^\natural(r)}{\partial r}\,dr = \lim_{\varepsilon\to 0}\varphi^\natural(\varepsilon) = \varphi(0). \qquad \diamond$$

# 8 Distributions and Fourier transform

How to understand the solutions (if any) of the PDE on generalized quaternion algebras (coquaternion's, double-complex etc)? Having in mind the classical meaning of the notion of symbols, we cannot understand directly the complex "hyperbolicity". We have in mind for instance the double-complex Laplacians $\Delta_+$ and $\Delta_-$. An answer is obtained for $\Delta_+$ by Petar Popivanov.

All things in the next of this paragraph are due to Petar Popivanov.

First, we decomplexify the operator

$$\Delta_+ = \partial^/\partial z^2 + i\partial^/\partial w^2,$$

where $z = x_0 + ix_1$, $w = x_2 + ix_3$, $i \in \mathbf{C}$, and $\partial/\partial z = \frac{1}{2}(\partial/\partial x_0 - i\partial/\partial x_1)$, $x' = (x_0, x_1)$, $x'' = (x_2, x_3) \in \mathbf{R}^2$.

We recall that a distribution $\varepsilon(x) \in \mathcal{D}'$ is a fundamental solution for $\Delta_+$ if

$$\Delta_+(D)\varepsilon(x) = \delta(x), \ x \in \mathbf{R}^4, \ \text{where } \delta(x) \ \text{is the Dirac delta function}.$$

The characteristic set of the operator ( . ) is defined as follows

$$Char\,\Delta_+(D) = \{(x,\xi) \in T^*(\mathbf{R}^4), \ (x,\xi) \neq 0 : p(\xi) = 0\}.$$

The symbol of $\Delta_+$ is denoted by $p$, $p(\xi) = \dfrac{1}{4}[(i\xi_1 + \xi_2)^2 + i(i\xi_3 + \xi_4)^2] = -\dfrac{1}{4}[(\xi_1 - i\xi_2)^2 + i(\xi_3 - i\xi_4)^2 = (\overline{\lambda})^2 + i(\overline{\mu})^2$, $\lambda := \xi_1 + i\xi_2$, $\mu := \xi_3 + i\xi_4$. Thus it is obtained $-4p(\xi) = (\overline{\lambda})^2 + i(\overline{\mu})^2$. As a consequence it is fulfilled

$$p(\xi) = 0 \iff \lambda^2 - i\mu^2 = 0 \iff \lambda = \pm e^{i/4}\mu \Rightarrow \lambda\mu \neq 0, \ \text{as } \lambda \neq 0, \mu \neq 0.$$

This shows that $\Delta_+$ can be written in a real and complex forms:

$$Char\,\Delta_+ = \mathbf{R}^4 \times Char_\xi\Delta_+,$$

where $Char_\xi\Delta_+ = \{\xi \in \mathbf{R}^4 \setminus 0 : \xi_1^2 - \xi_2^2 + 2\xi_3\xi_4 = 0, \ \xi_3^2 - \xi_4^2 - 2\xi_1\xi_2 = 0\} = \{(\lambda, \mu) \in \mathbf{C}^2 : \lambda = \pm e^{i\pi/4}\mu\}$ and $\lambda \neq 0, \mu \neq 0$.

**Lemma.** $Char_\xi\Delta_+$ is a 2-dimensional conic submanifold of $\mathbf{R}^4$.

**Theorem (P.Popivanov).** *The following statements are valid for $\Delta_+$ :*

*1. $\Delta_+$ is neither $C^\infty$, nor analytic hypoelliptic.*

*2. $\varepsilon(x) = c/(\overline{z}^2 + \overline{w}^2) \in \mathcal{D}'(\mathbf{R}^4)$, $c = const \neq 0$, is the fundamental solution for $\Delta_+$.*

*3. $\Delta_+$ is partially $C^\infty$ hypoelliptic.*

*4. $\Delta_+$ is operator of principal type. More generally, one say that the linear operator $P$ is an operator of principal type if $p(\xi) = 0 \Rightarrow \nabla_\xi p(\xi) \neq 0$. It is to remark that if $w \in \mathcal{D}'(\Omega)$ and $Pw = f$ then $(Char\,P \cap WF(w)) \setminus WF(f)$ is invariant under the bicharacteristic foliation in $Char\,P \setminus WF(f)$. This means that the singularities $Pw = f$, $f \in C^\infty$, propagate along the bicharacteristics $\Gamma$ of $p$. Here as ordinary, $WF(f)$ denote the wave front of the Schwartz distribution $f$.*

*5. If $Pw \in C^\infty$ and $w \in \mathcal{D}'(\Omega)$, then for each $\varphi \in C_0^\infty$ and each natural number $N$ there exists a positive constant $C_N$, such that* $|\hat{\varphi w}(\xi)| \leq C_N/(1+|\xi''|)^N$ when $\xi_1 + i\xi_2/\xi_3 + i\xi_4 \in M^1$, and respectively,
$|\hat{\varphi w}(\xi)| \leq D_N/(1+|\xi'|)^N$ when $\xi_3 + i\xi_4/\xi_1 + i\xi_2 \in M^2$ .

For $g \in \mathcal{D}'$, $\hat{g}(\xi)$ stands for Fourier transform of the function and $WF(f)$ is the wave front of the Schwartz distribution of $f \in \mathcal{D}'$.

Proof of p.2.

As it is known (see [V]) $\mathcal{E}_1(u) \equiv \mathcal{E}_1(x_0, x_1) = \dfrac{1}{\pi \overline{z}} = \dfrac{1}{\pi(x_0 - ix_1)}$ is a fundamental solution of $\dfrac{\partial}{\partial z}$, i. e.

$$\mathcal{E}_1'(u) = \frac{\partial}{\partial u}\mathcal{E}_1(x_0, x_1) = \delta(x_0, x_1) = \delta(Rez, Imz).$$

Let $a \in \mathbf{C}^1$. Then $\dfrac{\partial}{\partial u}\mathcal{E}_1(u+a) = \mathcal{E}_1'(u+a) = \delta((x_0, x_1) + (Rea, Ima))$.

It is looking for the fundamental solution of the form

$$\mathcal{E}(x) = \mathcal{E}(z, w) = \mathcal{E}_1(z + e^{i\pi/4}w)\mathcal{E}_1(z - e^{i\pi/4}w).$$

Therefore,

$$\Delta_+\mathcal{E} = 4\mathcal{E}_1'(z + e^{i\pi/4}w)\mathcal{E}_1'(z - e^{i\pi/4}w) =$$
$$= 4\delta(x' + (Ree^{i\pi/4}w, Ime^{i\pi/4}w))\delta(x' - (Ree^{i\pi/4}w, Ime^{i\pi/4}w)).$$

Put $y_1 = Re(e^{i\pi/4}w)$, $y_2 = Im(e^{i\pi/4}w)$. Then the identity $\delta(x' + y')\delta(x' - y') = const.\delta(x', y')$ completes the proof. (See ([H], Chapter 8)).

In fact, let $\varphi(x', y') \in C_0^\infty(\mathbf{R}^4)$. Then the linear change $x' - y' = z, x' + y' = w$ leads to

$$\langle \delta(x' - y')\delta(x' + y'), \varphi(x', y')\rangle = const.\langle \delta(z) \otimes \delta(w),$$

$$\varphi\left(\frac{z+w}{2}, \frac{w-z}{2}\right)\rangle = const.\varphi(0,0) = const.\langle \delta(x', y'), \varphi(x', y')\rangle;$$

$const. \neq 0$, $y' = 0 \iff e^{i\pi/4}w = 0 \iff w = 0 \iff x_2 = x_3 = 0$.

For the full proof of the theorem see [P].

# 9  Hyperbolic 4-real harmonicity

Here we shall consider hyperbolic 4-real numbers $x = x_0 + jx_1 + j^2x_2 + j^3x_3$, $j^4 = +1$, and functions of hyperbolic 4-real coordinates

$$f(x) = \sum_{k=0}^{3} j^k f_k(x_0, x_1, x_2, x_3).$$

We suppose that functions $f_k$ are continuously differentiable at least two times.

According to [ADMS], for the hyperbolic holomorphic 4-real function we have two couples of hyperbolic $CR$ equations, namely

$$\frac{\partial f_0}{\partial x_0} = \frac{\partial f_2}{\partial x_2}, \qquad \frac{\partial f_0}{\partial x_2} = \frac{\partial f_2}{\partial x_0} \quad (\text{ Hyperbolic } \ CR(x_0, x_2) \ \text{ system})$$

$$, \frac{\partial f_1}{\partial x_1} = \frac{\partial f_3}{\partial x_3}, \qquad \frac{\partial f_1}{\partial x_3} = \frac{\partial f_3}{\partial x_1} \quad (\text{ Hyperbolic } \ CR(x_1, x_3) \ \text{ system}).$$

Differentiating we obtain the following system of second order patrial differential equations

$$(a) \qquad \frac{\partial^2 f_0}{\partial x_0^2} - \frac{\partial^2 f_0}{\partial x_2^2} = 0, \quad \frac{\partial^2 f_2}{\partial x_0^2} - \frac{\partial^2 f_2}{\partial x_2^2} = 0,$$

for the functions $f_0$ and $f_2$ and

$$(b) \qquad \frac{\partial^2 f_1}{\partial x_1^2} - \frac{\partial^2 f_1}{\partial x_3^2} = 0, \quad \frac{\partial^2 f_3}{\partial x_1^2} - \frac{\partial^2 f_3}{\partial x_3^2} = 0,$$

for the functions $f_1$ and $f_3$.

The functions $f_0$ and $f_2$ are called hyperbolic harmonic with respect to the system a) and, respectively, the functions $f_1$ and $f_3$ are called hyperbolic harmonic with respect to the system b). The function $f$ is hyperbolic holomorphic.

**Example.** The function $f_0(x_0, x_1, x_2, x_3) = x_0^2 + x_2^2$ is hyperbolic harmonic with respect to to a) and respectively $f_1(x_0, x_1, x_2, x_3) = x_1^2 + x_3^2$ is hyperbolic harmonic with respect to to b). $f(x_0 + jx_1 + j^2 x_2 + j^3 x_3) = = x_0^2 + x_2^2 + j(x_1^2 + x_3^2) + j^2(2x_0 x_2) + j^3(2x_1 x_3)$ is hyperbolic holomorphic.

## 10   Hyperbolic 4-real geometry

The dual representation of the elements of the algebras $\mathbf{R}(1, j, j \cdot j^3)$, with $j^4 = -1$ and with $j^4 = +1$, as a linear forms with matrix coefficients of one side, and as special matrices $4 \times 4$ of other side, inspires two different point of view for their geometry. In the first case some remarks was presented in [ADMS] from the point of linear algebra about linear forms $x_0 + jx_1 + j^2 x_2 + j^3 x_3$. The second case concerns the matrix representation points. Especially for hyperbolic algebra $\mathbf{R}(1, j, j \cdot j^3)$ with $j^4 = +1$, we recall that

$$(29) \qquad x = x_0 + jx_1 + j^2 x_2 + j^3 x_3 = \begin{pmatrix} x_0 & x_1 & x_2 & x_3 \\ x_3 & x_0 & x_1 & x_2 \\ x_2 & x_3 & x_0 & x_1 \\ x_1 & x_2 & x_3 & x_0 \end{pmatrix}, \ j^4 = +1.$$

We will denote by $\Delta_h(x)$ the determinant of the matrix of $x$ above.

**Lemma.** (L. Apostolova) *The following formula is true*

$$(30) \qquad \Delta_h(x) + 4(x_0 x_3 - x_1 x_2)(x_0 x_1 - x_2 x_3) = \left(x_0^2 - x_2^2\right)^2 - \left(x_1^2 - x_3^2\right)^2.$$

The proof is obtained by straightforward calculation. $\qquad\qquad \diamond$

**Remark.** The lemma above generalizes the analogous lemma for hyperbolic 2-real numbers (in fact, ordinary hyperbolic complex numbers)

$$x = x_0 + jx_1 = \begin{pmatrix} x_0 & x_1 \\ x_1 & x_0 \end{pmatrix}, \ j^2 = +1,$$

namely

$$(31) \qquad \det \begin{pmatrix} x_0 & x_1 \\ x_1 & x_0 \end{pmatrix} = x_0^2 - x_1^2.$$

The form in the right side of (30), $\left(x_0^2 - x_2^2\right)^2 - \left(x_1^2 - x_3^2\right)^2$ is not a positive definite form like the analogous form in the right side of the equality (31). Evidently from (30) it follows that

$$-\left(x_1^2 - x_3^2\right)^2 \le \Delta_h(x) + 4(x_0 x_3 - x_1 x_2)(x_0 x_1 - x_2 x_3) \le \left(x_0^2 - x_2^2\right)^2$$

and the following estimation for the determinant $\Delta_h(x)$ hold

$$-\left(x_1^2 - x_3^2\right)^2 - 4(x_0 x_3 - x_1 x_2)(x_0 x_1 - x_2 x_3) \le \Delta_h(x) \le$$

$$\le \left(x_0^2 - x_2^2\right)^2 - 4(x_0 x_3 - x_1 x_2)(x_0 x_1 - x_2 x_3).$$

It is to have in mind also the inequality

$$|(x_0^2 - x_2^2)^2 - (x_1^2 - x_3^2)^2| \le (x_0^2 - x_2^2)^2 + (x_1^2 - x_3^2)^2.$$

The determinant $\Delta_h(x)$ is a polynomial of degree 4. It is called characteristic polynomial of the algebra $\mathbf{R}(1, j, j \cdot j^3)$, $j^4 = +1$. We have received, in fact, an estimation for the mentioned characteristic polynomial.

The polynomial $P_h(x)$ determines the hyperbolic 4-real geometry.

## Appendix A   (Stancho Dimiev)

The notion of double-complex number was introduced in [D] as a special alternative of the notion of bicomplex numbers. Here we expose a motivation. The double-complex numbers arise naturally from the 4-real numbers as follows: taking a 4-real number $x = x_0 + jx_1 + j^2 x_2 + j^3 x_3$ one set $x_0 + j^2 x_2 = z_0$ and $x_1 + j^2 x_3 = z_1$. So $x = z_0 + jz_1$. As $(j^2)^2 = -1$ we identify $j^2$ with the complex number $i$, $i \in \mathbf{C}$. In such a way the square of $j$ in the expression $x = z_0 + jz_1$ became equal to $i$, $i \in \mathbf{C}$. The mentioned expression $z_0 + jz_1$ determine a double-complex number.

Analogously, the notion of 4-complex number arise from the 8-real numbers $x = \sum_{k=0}^{3} j^k x_k$. Clearly, $x = (x_0 + j^4 x_4) + j(x_1 + j^5 x_5) + j^2(x_2 + j^4 x_6) + j^3(x_3 + j^4 x_7)$. Setting $Z_0 := x_0 + j^4 x_4$, $Z_1 := x_1 + j^4 x_5$, $Z_2 := x_2 + j^4 x_6$, $Z_3 := x_3 + j^4 x_7$ we receive

$$x = Z_0 + jZ_1 + j^2 Z_2 + j^3 Z_3,$$

where $Z_k$, $k = 0, 1, 2, 3$, can be considered as a double-complex numbers (as $(j^4)^2 = -1$ we can identify $j^4$ with the complex number $i$).

The following scheme of generalizations appear

I. Complex numbers $z, w, \dots$ .

II. Double-complex number $z + jw, \ldots$, where $z$ and $w$ are complex numbers and $j^2 = i$.

III. Four-complex (4-complex) numbers $Z_0 + jZ_1 + j^2Z_2 + j^3Z_3$, $j^4 = = i$ as $j^8 = -1$, where $Z_0, Z_1, Z_2, Z_3$ can be considered as double-complex numbers.

IV Eight-complex (8-complex) numbers $W_0 + jW_1 + \ldots + j^7W_7$, $j^8 = i$, as $j^{16} = -1$ where $W_k$ can be considered as a four-complex numbers.

etc.

The Cauchy formula is easily formulated for the double-complex numbers and by induction for the four-complex numbers etc.

The Weierstrass Preparation Theorem holds for double-complex numbers and so on.

We are interested to have a type of Weierstrass Preparation Theorem for pseudodifferential operators on hypercomplex algebras.

Analogously, we define hyperbolic 4-real matrix numbers with the help of matrix $j = \mathcal{H}_{1,3}$ and its degrees $j^2$, $j^3$ and $j^4 = (+1)\mathrm{I}_4$. Now we set $j^2 = \mathbf{j}$, from where $\mathbf{j}^2 = j^4 = (+1)\mathrm{I}_4$. Identifying $\mathrm{I}_4$ with 1, we set $x_0 + j^2x_2 = = x_0 + \mathbf{j}x_2 = Z$, $x_1 + j^2x_3 = x_1 + \mathbf{j}x_3 = W$, we obtain the expression

$$Z + jW = (x_0 + \mathbf{j}x_2) + j(x_1 + \mathbf{j}x_3),$$

that can be view as a hyperbolic double-complex number as $Z$ and $W$ can be considered as hyperbolic complex numbers and $j^2 = \mathbf{j}$.

Analogous argument as for 4-complex numbers $z_0 + jz_1 + j^2z_2 + j^3z_3$, $j^4 = \mathbf{j}$, $\mathbf{j}^2 = 1$, $j^8 = -1$, can be used for construction of hyperbolic four-complex numbers etc.

In general, it is important for us to describe precisely the geometry of hyperbolic 4-real numbers and this one of hyperbolic double-complex numbers.

Let us remark that the ordinary hyperbolic complex numbers $(t + jx, j^2 = 1)$ possess a natural Minkowski metric

$$d\tau^2 = dt^2 - dx^2$$

(in our matrix notation $\begin{pmatrix} x_0 & x_1 \\ x_1 & x_0 \end{pmatrix}$, we have $d\tau^2 = dx_0^2 - dx_1^2$), but in the case of hyperbolic 4-real numbers there are two natural Minkowski metric

$$d\tau_0^2 = dx_0^2 - dx_2^2, \quad \text{and} \quad d\tau_1^2 = dx_1^2 - dx_3^2$$

(to see the paragraph 9 above).

In general, it is important to describe precisely the geometry of hyperbolic 4-real variables, having in mind also hyperbolic 4-real manifolds, and the corresponding one of hyperbolic double-complex variables (differential forms, manifolds, etc.).

A new aspect appears naturally. From one side: the (elliptic) 4-real variables with Euclidean metric and the corresponding double-complex variables with an Hermitean double-complex metric.

From other side: the hyperbolic 4-real variables with Minkowski signature and related Minkowski metric, and the corresponding geometry of the hyperbolic double-complex variables (coordinate).

This means to introduce the notion of almost double-complex structure on two-dimensional double-complex manifolds and, respectively, the notion of hyperbolic double-complex structure.

It is interesting to clarify the role of the Weierstrass preparation theorem for pseudo-differential operators in the above described geometries.

## Appendix B   (L. N. Apostolova)

*1. Euler formula for the double-complex exponential function* (see [A2])

The exponential function $e^{z+jw}$ is well defined function for double-complex numbers by the absolutely convergent power series $e^{\alpha} = \sum_{n=0}^{\infty} \alpha^n/n!$, $\alpha \in \mathbf{C}(1,j)$.

Some important properties of $e^{\alpha}$ are the following ones:

- if we define a formal derivative with respect to the double-complex variable $\alpha = z + jw$ to be the following one: $\dfrac{\partial}{\partial \alpha} := \dfrac{1}{2}\left(\dfrac{\partial}{\partial z} - ij\dfrac{\partial}{\partial w}\right)$, where $\dfrac{\partial}{\partial z}$ and $\dfrac{\partial}{\partial w}$ are the formal complex derivative with respect to complex variables $z$ and $w$, then the first derivative of the exponential double-complex function $e^{\alpha}$ is the same function $e^{\alpha}$;

- the formula $e^{a+b} = e^a e^b$ is valid;

- the values of the exponential function are double-complex numbers, which are not divisors of the zero;

- it is valid a formula analogous to the Euler formula for complex numbers $e^{i\varphi} = \cos\varphi + i\sin\varphi$ ($\varphi$-real number, $i$ - the imaginary unit).

We shall derive formulas for the functions $C(w)$, $S(w) : \mathbf{C}(1,j) \to \mathbf{C}$, where $e^{z+jw} = e^z(C(w) + jS(w))$, $(z,w) \in \mathbf{C}(1,j)$. From the series development of the function $e^{jw}$ we obtain

$$e^{jw} = \sum_{n=0}^{\infty} \frac{(jw)^n}{n!} = \sum_{k=0}^{\infty} \frac{(jw)^{2k}}{(2k)!} + j\sum_{k=0}^{\infty} \frac{j^{2k}(w)^{2k+1}}{(2k+1)!} =$$

$$= \sum_{k=0}^{\infty} \frac{i^k(w)^{2k}}{(2k)!} + j\sum_{k=0}^{\infty} \frac{i^k(w)^{2k+1}}{(2k+1)!} =: C(w) + jS(w).$$

The functions $C(w)$ and $S(w)$ are entire functions of one complex variable, called an even and an odd part of the function $e^{iw}$. The function $C(w)$ is an even function and $S(w)$ is an odd function of the variable $w$. It is fulfilled $C(w) = (e^{jw} + e^{-jw})/2$ and $S(w) = (e^{jw} - e^{-jw})/2$.

To express these functions with the complex trigonometric functions $\cos w$ and $\sin w$ let us check the compose function

$$C\left(\frac{1+i}{\sqrt{2}}w\right) = \sum_{k=0}^{\infty} i^k \left(\frac{1+i}{\sqrt{2}}\right)^{2k} w^{2k} = \sum_{k=0}^{\infty} \frac{(-1)^k w^{2k}}{(2k)!} = \cos w.$$

Analogously we obtain

$$S\left(\frac{1+i}{\sqrt{2}}w\right) = \sum_{k=0}^{\infty} i^k \left(\frac{1+i}{\sqrt{2}}\right)^{2k+1} w^{2k+1} = \frac{1+i}{\sqrt{2}} \sum_{k=0}^{\infty} \frac{(-1)^k w^{2k+1}}{(2k+1)!} = \frac{1+i}{\sqrt{2}}\sin w.$$

So the following formulas arise

$$C\left(\frac{1+i}{\sqrt{2}}w\right) = \cos w, \quad S\left(\frac{1+i}{\sqrt{2}}w\right) = \frac{1+i}{\sqrt{2}}\sin w \quad \text{and} \tag{1}$$

$$C(w) = \cos\frac{1-i}{\sqrt{2}}w, \quad S(w) = \frac{1+i}{\sqrt{2}}\sin\frac{1-i}{\sqrt{2}}w.$$

The functions $C(w)$ and $S(w)$ are periodic functions with period equal to $\sqrt{2}(1-i)\pi$ and the identity $C^2(w) - iS^2(w) = 0$ holds.

The formulas for the derivatives of $C(w)$ and $S(w)$ hold

$$C'(w) = -\left(\sin\frac{1-i}{\sqrt{2}}w\right)\frac{1-i}{\sqrt{2}} = iS(w) \quad \text{and}$$

$$S'(w) = \frac{1+i}{\sqrt{2}}\left(\cos\frac{1-i}{\sqrt{2}}w\right)\frac{1-i}{\sqrt{2}} = \cos\frac{1-i}{\sqrt{2}}w = C(w).$$

From this two formulas it is easy to obtain the fundamental equality

$$\frac{d}{dw}e^{jw} = je^{jw}.$$

Let us prove that the double-complex exponential function is a holomorphic double-complex function (see also [A2]). To do this we check

$$\frac{1}{2}\left(\frac{\partial e^{z+jw}}{\partial z} + ij\frac{\partial e^{z+jw}}{\partial w}\right) = e^{z+jw} + ijje^{z+jw} = e^{z+jw} - e^{z+jw} = 0.$$

Moreover, the function $e^{z+jw}$ is a harmonic double-complex function. Indeed,

$$\Delta_+ e^{z+jw} = \frac{\partial^2 e^{z+jw}}{\partial z^2} + i\frac{\partial^2 e^{z+jw}}{\partial w^2} = e^{z+jw} + ij^2 e^{z+jw} = 0.$$

*2. Euler formula for the hyperbolic double-complex exponential function* (see [A3])

The exponential function $e^{z+jw}$ is well defined function of hyperbolic double complex numbers $\alpha = z + jw = x_0 + jx_1 + j^2x_2 + j^3x_4$, $x_k \in \mathbf{R}$, $k = 0, 1, 2, 3$ with the absolutely convergent power series $e^\alpha = \sum_{n=0}^{\infty}\alpha^n/n!$, $\alpha \in D\tilde{\mathbf{C}}_2)$, $j^4 = +1$. We suppose here that the algebra $D\tilde{\mathbf{C}}_2$ is equipped with the scalar product $\langle x_0 + jx_1 + j^2x_2 + j^3x_3, y_0 + jy_1 + j^2y_2 + j^3y_3 \rangle = |x_0y_0 - x_1y_1 + x_2y_2 - x_3y_3|$.

Some important properties of $e^\alpha$ are the following ones:

- the first derivative of the exponential hyperbolic double-complex function with respect to it hyperbolic double-complex variable is the same function;

- the formula $e^{a+b} = e^a e^b$ is valid for $a, b$ - hyperbolic double-complex variables;

- the values of the exponential function are hyperbolic double-complex numbers, which are not divisors of the zero;

- it arise a formula for hyperbolic double-complex numbers, analogous to the Euler formula for complex numbers $e^{i\varphi} = \cos\varphi + i\sin\varphi$ ($\varphi$-real number, $i$ - the imaginary unit).

Let us recall the series development of the real functions $e^x$, $\cos x$, $\sin x$, $\cosh x$ and $\sinh x$. It is well known that

$$e^x = \sum_{k=0}^{\infty} \frac{x^k}{k!}, \quad \cos x = \sum_{k=0}^{\infty} \frac{(-1)^k x^{2k}}{(2k)!}, \quad \sin x = \sum_{k=0}^{\infty} \frac{(-1)^k x^{2k+1}}{(2k+1)!},$$

$$\cosh x = \sum_{k=0}^{\infty} \frac{x^{2k}}{(2k)!}, \quad \sinh x = \sum_{k=0}^{\infty} \frac{x^{2k+1}}{(2k+1)!},$$

from where follows

$$\cosh x + \cos x = \sum_{k=0}^{\infty} \frac{x^{4k}}{(4k)!}, \quad \sinh x + \sin x = \sum_{k=0}^{\infty} \frac{x^{4k+1}}{(4k+1)!},$$

$$\cosh x - \cos x = \sum_{k=0}^{\infty} \frac{x^{4k+2}}{(4k+2)!}, \quad \sinh x - \sin x = \sum_{k=0}^{\infty} \frac{x^{4k+3}}{(4k+3)!}.$$

By the developments of the series of the hyperbolic double-complex function, using the polynomial formula we obtain

$$e^{\alpha} := \sum_{n=0}^{\infty} (z+jw)^n / n! = \sum_{k_0,k_1,k_2,k_3 \geq 0} \frac{1}{k_0! k_1! k_2! k_3!} x_0^{k_0} (jx_1)^{k_1} (j^2 x_2)^{k_2} (j^3 x_3)^{k_3} =$$

$$= e^{x_0} \sum_{k=0}^{\infty} \left( \frac{1}{(4k)!} x_1^{4k} + \frac{j}{(4k+1)!} x_1^{4k+1} + \frac{j^2}{(4k+2)!} x_1^{4k+2} + \frac{j^3}{(4k+3)!} x_1^{4k+3} \right) \times$$

$$\times \sum_{k=0}^{\infty} \left( \frac{1}{(2k)!} x_2^{2k} + j^2 \frac{1}{(2k+1)!} x_2^{2k+1} \right) \times$$

$$\times \sum_{k \geq 0} \left( \frac{j^{12k}}{(4k)!} x_3^{4k} + \frac{j^{12k+3}}{(4k+1)!} x_3^{4k+1} + \frac{j^{12k+6}}{(4k+2)!} x_3^{4k+2} + \frac{j^{12k+9}}{(4k+3)!} x_3^{4k+3} \right) =$$

$$= e^{x_0} (\cosh x_2 + j^2 \sinh x_2) \times$$

$$\times ((1+j^2) \cosh x_1 + (j+j^3) \sinh x_1 + (1-j^2) \cos x_1 + (j-j^3) \sin x_1) \times$$

$$\times ((1+j^2) \cosh x_3 + (j^3+j) \sinh x_3 + (1-j^2) \cos x_3 + (j^3-j) \sin x_3).$$

Let us recall that the hyperbolic double-complex numbers $1+j^2$ and $1-j^2$ are divisors of zero and it is true that $(1+j^2)(1-j^2) = 0$.

Let us check the action of the "hyperbolic" double-complex Laplace operator on the hyperbolic double-complex exponential function. We obtain

$$\Delta_h e^{\alpha} = \frac{\partial^2 e^{\zeta+j\eta}}{\partial \zeta^2} - \mathbf{j} \frac{\partial^2 e^{\zeta+j\eta}}{\partial \eta^2} = e^{\zeta+j\eta} - j^2 \mathbf{j} \, e^{\zeta+j\eta} = 0,$$

i.e. the even and the odd parts of the exponential hyperbolic double-complex function $e^{\alpha}$ are harmonic hyperbolic complex-valued functions.

We will prove also that the function $e^{\alpha}$ is a holomorphic hyperbolic double-complex function. To do this we check

$$\left( \frac{\partial e^{\zeta+j\eta}}{\partial \zeta} - \frac{1}{j} \frac{\partial e^{\zeta+j\eta}}{\partial \eta} \right) = e^{\zeta+j\eta} - j\frac{1}{j} e^{\zeta+j\eta} = 0.$$